\documentclass[12pt]{article}
\usepackage{amsmath}
\usepackage{amsfonts}
\usepackage{amsthm}
\usepackage{graphicx}
\usepackage{bbm}
\usepackage{colortbl}
\theoremstyle{plain}

\newtheorem{proposition}{Proposition}

\theoremstyle{definition}
\newtheorem{definition}{Definition}

\newtheorem{example}{Example}

\newcommand{\bs}{\boldsymbol}

\title{Equilibrium strategies in a multiregional transboundary pollution differential game with spatially distributed controls\thanks{
This research is partially supported by MEC under projects
 MTM2016-78995-P (AEI)  (first author) and   ECO2017-82227-P (AEI) (third author), and by Junta de Castilla y Le\'on VA024P17 and  VA105G18 co-financed by FEDER funds (EU). }}

\author{Javier de Frutos \and Paula L\'{o}pez \and Guiomar Mart\'{\i}n-Herr\'{a}n}
\date{\small
IMUVA, Universidad de Valladolid. Spain. }

\begin{document}
\maketitle
\abstract{We analyze a differential game with spatially distributed controls to study a multiregional transboundary pollution problem. The dynamics of the state variable (pollution stock) is defined by a two dimensional parabolic partial differential equation. The control variables (emissions) are spatially distributed variables. The model allows for a, possibly large, number of agents with predetermined geographical relationships. For a special functional form previously used in the literature of transboundary pollution dynamic games we analytically characterize the feedback Nash equilibrium. We show that at the equilibrium both the level and the location of emissions of each region depend on the particular geographical relationship among agents. We prove that, even in a simplified model, the geographical considerations can modify the players' optimal strategies and therefore, the spatial aspects of the model should not be overlooked.

\vspace {0.5cm}
\noindent
{\bf Keywords}: Transboundary Pollution; Spatial Dynamics; Spatially Distributed Controls;  Differential Games; Parabolic Differential Equations.

\vspace {0.25cm}
\noindent

\setlength{\parskip}{\medskipamount}

\section{Introduction}

Most of the literature on dynamic models settled for the analysis  of different economic and environmental problems takes into account the temporal aspect but disregards the spatial aspect. The addition of the spatial dimension although enriches the model and consequently its prescriptions comes at the cost of making the analysis much more difficult.  Over the last decade the spatial dimension has been introduced in different economic and environmental contexts, such as, for example, allocation of economic activity or technological diffusion (Brito (2004), Camacho et al. (2008), Brock \& Xepapadeas (2008a), Boucekkine et al. (2009, 2013a, 2013b, 2019a),   Desmet \& Rossi-Hansberg (2010), Brock et al. (2014a) and Fabbri (2016)) or environmental and climate problems (Brock \& Xepapadeas (2008b, 2010), Brock et al. (2014b), Camacho \& P\'erez-Barahona (2015), Xepapadeas (2010), Anita et al. (2013, 2015),   Desmet \& Rossi-Hansberg (2015), La Torre et al. (2015), Augeraud-V\'eron et al. (2017, 2019a, 2019b), Boucekkine et al. (2019b), and De Frutos \& Mart\'\i n-Herr\'an (2019a, 2019b)). Except De Frutos \& Mart\'\i n-Herr\'an (2019a, 2019b) all the previously cited papers study finite or infinite time horizon optimal control problems  extended to infinite dimensional state space. All these papers analyze problems where there is only one decision-maker.

One of the main differences between De Frutos \& Mart\'\i n-Herr\'an (2019a, 2019b) and the other contributions to this literature is the introduction of strategic interactions between the decision-makers. The addition of the strategic behavior of the agents implies a  methodological change, moving from  an optimal control  problem to a dynamic game. More specifically, De Frutos \& Mart\'\i n-Herr\'an (2019a, 2019b) study an intertemporal transboundary pollution dynamic game where there is a continuum of spatial sites and the pollution stock diffuses over these sites. In both papers the original $J$-player (region or country) model is formulated   in continuous space and continuous time with two spatial dimensions and one temporal dimension, the  spatio-temporal evolution of the stock of a pollutant is described by a diffusion partial differential equation (PDE) and general boundary conditions are assumed. In De Frutos \& Mart\'\i n-Herr\'an (2019a) each player decides the emission level so as to maximize  the present value of benefits net of environmental damages due to the concentration of a pollutant over his spatial domain, taking into account the PDE that describes the spatio-temporal evolution of the stock of a pollutant. In De Frutos \& Mart\'\i n-Herr\'an (2019b)  in addition to choosing the level of emissions, each player also decides the investment in clean technology, taking into account the temporal evolution of the stock of clean technology and that the greater this stock the lower the emission-output ratio.

Both papers follow the same spatial discretization approach to characterize the equilibrium emission strategies in a multiregional setting with spatial effects. In the space-discretized model there are $J$ state variables,  the average pollution in each one of the $J$ regions or countries, and the temporal dynamics of these variables is described by a system of $J$ ordinary differential equations. Gra{\ss} \& Uecker (2017) uses a similar spatial discretization approach to analyze spatially distributed optimal control problems. In De Frutos \& Mart\'\i n-Herr\'an (2019a) the feedback Nash equilibrium of the two-player space-discretized model is analytically characterized, while for more than two players, the model is numerically solved by adapting a numerical algorithm presented in  De Frutos \& Mart\'\i n-Herr\'an (2015). The linear-quadratic specification of the model is inspired in the literature of transboundary pollution dynamic games (J{\o}rgensen et al. (2010) surveyed this literature), specifically in the seminal papers by Dockner \& Long (1993) and Van der Ploeg and De Zeeuw (1992). The numerical examples show that once the spatial dimension and the strategic behavior of the decision-makers are introduced the environmental policies greatly differ from those characterized ignoring either the spatial dimension or the strategic interactions among the agents. In De Frutos \& Mart\'\i n-Herr\'an (2019b)  particular functional forms  borrowed from J{\o}rgensen \& Zaccour (2001) are considered in such a way that the dynamic game belongs to the class of linear-state differential games. For this formulation the space-discretized model is exactly solved and is proved to be a clear generalization of the model which ignores the spatial transport phenomena. The equilibrium environmental policy of the spatial model  coincides with the equilibrium policy of the non-spatial model when the diffusion parameter, that describes how pollution diffuses among regions, tends to infinity.

In the present paper we revisit the study of the equilibrium emission strategies in a multiregional dynamic game  in a spatial context. The main objective of the paper is to characterize the feedback Nash equilibrium emission strategies of the original $J$-player (region or country) model  formulated   in continuous space and continuous time. We depart from the spatial discretization approach presented in De Frutos \& Mart\'\i n-Herr\'an (2019a, 2019b) and for a linear-state specification inspired in J{\o}rgensen \& Zaccour (2001) and De Frutos \& Mart\'\i n-Herr\'an (2019b) we  explicitly solve a functional Hamilton-Jacobi-Bellman system in an infinite dimensional Hilbert space. Thanks to the linear-state framework, guessing an affine form of the players' value functions we can explicitly
compute these functions and as a consequence also the feedback Nash equilibrium emission strategies. The resolution method involves solving several elliptical problems, one for each player.  Although the equilibrium strategies are constant in time as a direct consequence of the linear-state structure of the dynamic game, interestingly, these strategies are not constant in space. Solving the original model formulation enables to strategically determine the average total  emission in each region (as was the case in the  space-discretized formulation of the model in De Frutos \& Mart\'\i n-Herr\'an (2019a, 2019b)), but also to characterize  the  particular point in the space where each region is emitting  pollutants. Using the new approach proposed in this paper we can analyze the optimal intraregional distribution of emissions of the pollutant, a question that has been previously neglected.  Furthermore, the new approach allows the study not only of the diffusive effect, but also the case where advection is important and particles are transported due to some external convective field, as for example, wind or water flow.

Our analytical results show that the geographical aspects are essential ingredients when determining the equilibrium emission strategies. Through several examples with different geographical configurations we show  that the type and behavior of the neighbors of one region have an influence in its own equilibrium environmental policies, not only by determining the optimal amount of emissions but also establishing the optimal spatial location. The results corroborate those obtained in the analysis of the space-discretized model which allows us to conclude that this simplified model correctly captures the spatial essence of the setting. However, the original model analyzed in the present paper is richer which in turn allows us to have an accurate view of how the spatial ingredients affect the equilibrium emission strategies and hence the pollution stock.

The rest of the paper is organized as follows. Section 2 presents the multiregional transboundary pollution differential game with spatially distributed controls and introduces some technical hypotheses and definitions. Section 3 analytically studies the Hamilton-Jacobi-Bellman system of equations and derives the feedback Nash equilibrium of the differential game. Section 4 presents some examples to highlight the properties of the Nash equilibria. The paper finishes with some concluding remarks.

\section{The model}\label{TheModel}

Let us denote by $\Omega$ a bounded planar domain endowed with a partition  $\Omega_j$, $j=1,\dots, J$, such that
\begin{equation*}\label{partition}
\overline{\Omega}=\bigcup_{j=1}^J\overline{\Omega}_j, \quad \Omega_i\cap\Omega_j=\emptyset, \quad i\ne  j.
\end{equation*}
where  $\overline{\Omega}$ is the closure of $\Omega$. We denote by $\partial_{ij}$ the common boundary between subdomains $\Omega_i$ and $\Omega_j$, that is $\partial_{ij}:=\partial\Omega_i\cap\partial\Omega_j=\overline{\Omega}_i\cap\overline{\Omega}_j,\; i\ne  j$.

The model is a $J$-player differential game in which
the control variable of player $i$ is the pollutant emissions in region $\Omega_i$.  The game is played non-cooperatively. In what follows we identify player $i$ with the region $\Omega_i$ and use the word country to distinguish $\Omega_i$ from the whole region $\Omega$. We assume that each of the countries, $\Omega_i$, $i=1,\dots, J$,  can exclusively emit pollution in its own territory. The pollutant emissions in $\Omega_i$ are represented by a function $u_i :\Omega_i\times [0,+\infty)\rightarrow \mathbb{R}_+$,  $i=1,\dots, J$.

The objective of player $i$, $i=1,\dots, J$, is to  maximize the following functional with respect to $u_i$,
\begin{equation}\label{objective}
J_i(u_1,\dots,u_J,P_0)=\int_0^{+\infty}\int_{\Omega_i} e^{-\rho t}\Bigl(\log(u_i)-\varphi_i P\Bigr)\,d{\bs{x}}\,dt
\end{equation}
subject to
\begin{equation}\label{state_equation}
\begin{aligned}
&\frac{\partial P}{\partial t}=\nabla\cdot(k\nabla P)+\bs{b}\cdot\nabla P - c P + F(u_1,\dots,u_J),\quad \bs{x}\in\Omega,\\
&P(\bs{x},0)=P_0(\bs{x}),\quad\bs{x}\in\Omega,\\
&\alpha P(\bs{x},t)+k\nabla P\cdot\bs{n}=\alpha P_b(\bs{x},t),\quad \bs{x}\in\partial\Omega.
\end{aligned}
\end{equation}
Here $P$ denotes the state variable (pollution stock), $F$ is assumed to be a real smooth function of  its arguments, $P_0$ and $P_b$ are known functions representing respectively, the initial distribution of pollution and the external pollution  that enters into $\Omega$ through its boundary. In what follow we assume  without loss of generality that function $P_b$ is identically zero. The diffusion coefficient $k(\bs{x})$  is assumed to be a smooth function of the spatial variables satisfying $k_1\le k(\bs{x})\le k_2$ for some positive constants $k_1$ and $k_2$. The coefficient $c(\bs{x}) $ represents the natural decay of pollution and it is supposed to be a non-negative smooth function. Function $\bs{b}(\bs{x})$ represents an external convective field. We assume that $\bs{b}(\bs{x})$ is smooth and divergence free $\nabla\cdot\bs{b}=0$. Finally, $\alpha$ is a positive real number and $m_i$ denotes the area of  $\Omega_i$.

In this paper, we assume that $F(u_1,\dots,u_J)=\sum_{j=1}^J u_j\mathbbm{1}_{\Omega_j}$ where
$\mathbbm{1}_{\Omega_j}$ is the characteristic function of the set $\Omega_j$, $j=1,\dots, J$. We  are implicitly considering that function $u_i$ can be arbitrarily extended outside of the domain $\Omega_i$.

Let us consider $\mathbb{X}=L^2(\Omega)$  the Hilbert space of square integrable functions defined over domain $\Omega$. As usual we identify $\mathbb{X}$ with its dual $\mathbb{X}^*$. Here and in the rest of the paper, we denote by $H^s(\Omega)$ the Sobolev space of functions with $s$ derivatives, in the distributional sense, in $L^2(\Omega)$.

We define the set $\mathcal{U}_i$ of admissible controls for player $i$, $i=1,\dots, J$, as the set of functions $u_i$ defined in $\Omega_i\times\mathbb{R}_+\rightarrow\mathbb{R}_+$. Although it is not essential here, we also assume  that  functions $u_i\in \mathcal{U}_i$ are bounded for almost all $\bs{x}\in\Omega$ and $t\ge 0$.

 The hypotheses above guarantee that given an initial condition $P_0\in\mathbb{X}$, the state equation (\ref{state_equation}) has a unique weak solution for each choice of the  controls $u_i\in\mathcal{U}_i$, $i=1,\dots,J$, see Barbu (1993), Li \& Yong (1995), Tr\"{o}ltzsch (2009).

With this model we adopt the simplest version of the economic and environmental model that still presents two important features  allowing us to answer our main research question. First, the strategic behaviour of the players, emissions by one player affects the environment of all; and second, the spatial aspect that allows us to show that at the equilibrium both the level and the location of emissions of each region depend on the particular geographical relationship among agents.

We are interested in stationary Markov-perfect Nash equilibria (MPNE) of the game. Then, we look for controls of the form  $u_i(\bs{x},t)=\Lambda_i(P(\bs{x},t))$, $i=1,\dots, J$. Here the strategies $\Lambda_i$ are functionals $\Lambda_i:\mathbb{X}\rightarrow\mathcal{U}_i$ such that the controlled dynamics
\begin{equation}\label{dynamics_controlled}
\begin{aligned}
&\frac{\partial P}{\partial t}=\nabla\cdot(k\nabla P)+\bs{b}\cdot\nabla P - c P + \sum_{j=1}^J\Lambda_j(P)\mathbbm{1}_{\Omega_j},\quad \bs{x}\in\Omega,\\
&P(\bs{x},\tau)=P_\tau(\bs{x}),\quad\bs{x}\in\Omega,\\
&\alpha P(\bs{x},t)+k\nabla P\cdot\bs{n}=\alpha P_b(\bs{x},t),\quad \bs{x}\in\partial\Omega,
\end{aligned}
\end{equation}
has a unique solution defined in $[\tau,\infty)$ for every $\tau\ge 0$ and $P_\tau\in\mathbb{X}$.
\begin{definition}
A vector $\bs{\Lambda}^*=[{\Lambda}_1^*,\dots,{\Lambda}_J^*]$ of admissible strategies is a Markov Perfect Nash Equilibrium if $J_i(\bs{u}^*,P_0)\ge J_i([u_i,\bs{u}^*_{-i}],P_0)$, for all $u_i=\Lambda_i(P)$ with  $\Lambda_i$ an admissible strategy, $i=1,\dots, J$. Here,
$\bs{u}^*=[u_1^*,\dots,u_J^*]$, $u_j^*=\Lambda_j^*(P^*)$ and  $P^*$ is the solution of (\ref{dynamics_controlled}) with $\Lambda_i=\Lambda_i^*$, $i=1,\dots, J$. We use $[u_i,\bs{u}_{-i}^*]$ to denote $[u_1^*,\dots,u_i,\dots,u_J^*]$.

Given a stationary  MPNE $\bs{\Lambda}^*=[\Lambda_1^*,\dots,\Lambda_J^*]$, $V_i(P) = J_i(\bs{u}^*,P)$ is called the value function of Player $i$.
\end{definition}

The objective of the next section is to characterize the MPNE of the differential game (\ref{objective})-(\ref{state_equation}) through the study of the value function.

\section{The Hamilton-Jacobi-Bellman equation}

 For simplicity in the exposition we assume in this section that $\bs{b}\cdot\bs{n}=0$ in $\partial\Omega$.   We define the linear operator $\bs{A}:D(\bs{A})\rightarrow X$ by
$$\bs{A}P=\nabla\cdot(k\nabla P)+\bs{b}\cdot\nabla P - c P,\quad \forall P\in \mathbb{X},$$
where $D(\bs{A})=\{P\in H^2(\Omega)|\alpha P(\bs{x})+ k\nabla P\cdot\bs{n}=0 \},$
is the domain  of $\bs{A}$.
The linear operator $\bs{A}$ is continuous from $D(\bs{A})$ in $\mathbb{X}$ and  is the infinitesimal generator of a contraction semigroup $e^{At}$ in $\mathbb{X}$, see Li \& Yong (1995). In what follows $\nabla W(P)$ denotes the Fr\'{e}chet derivative of a functional $W$ and  $\langle \cdot ,\cdot\rangle$ denotes the scalar product in $\mathbb{X}$.

The following proposition is a consequence of Proposition 1.2 (Ch. 6, p.225) in Li \& Yong (1995), see also
Ba\c{s}ar \& Olsder (1999), Haurie et al. (2012).

\begin{proposition}\label{proposition_HJB}
Let $\bs{\Lambda}^*=[\Lambda_1^*,\dots,\Lambda_J^*]$ be a MPNE. Let us assume that $V^i(P)$ is of class $\mathcal{C}^1(\mathbb{X})$, $i=1,\dots, J$. The value functions $V^i$, $i=1,\dots, J$ satisfy the functional Hamilton-Jacobi-Bellman system
\begin{equation}\label{HJB}
\rho V^{i}(P)=
\sup_{u_i}\Bigl\{\mathcal{G}^{i}(P,u_i)+\big\langle \bs{A}P+u_i\mathbbm{1}_{\Omega_i}+\sum_{j\neq i} \Lambda^*_j(P)\mathbbm{1}_{\Omega_j},\nabla V^{i}(P)\big\rangle\Bigr\},\quad i=1,\dots, J,
\end{equation}
where
$$\mathcal{G}^{i}(P,u_i)=\int_{\Omega_i} \Bigl(\log(u_i)-\varphi_i P\Bigr)\,d{\bs{x}}.$$
Furthermore, $\Lambda_i^*(P)$, $i=1,\dots, J$, is a maximizer of the right hand side of (\ref{HJB}).
\end{proposition}

The regularity hypothesis in Proposition \ref{proposition_HJB} is very demanding. It is well known that even in the finite dimensional case Hamilton-Jacobi-Bellman equations can fail to have enough regularity and one has to resort to weaker concepts as viscosity solutions, see, for example, Barbu (1993), Cannarsa \& Da Prato (1990), Li \& Yon (1995). However,  Proposition \ref{proposition_HJB} is enough for the model we study in this paper.

In this paper, for simplicity, we use the strong optimality concept, Dockner et al. (2000), so that  transversality conditions of the form
\begin{equation}\label{transversality}
\lim_{t\mapsto\infty}e^{-\rho t}V_i(P^*(\bs{x},t))=0,\quad i=1,\dots,J.
\end{equation}
are necessary. The equilibrium we are going to compute explicitly leads to controlled dynamics that possess a (unique) stable stationary state, so that the transversality conditions will be automatically satisfied.

We remark that, in general, Hamilton-Jacobi-Bellman systems of the form (\ref{HJB}) with boundary conditions as (\ref{transversality})  have multiple solutions, see De Frutos \& Mart\'{\i}n-Herr\'{a}n (2018), which correspond with possible multiple MPNE. Because the differential game fits the linear-state class, we concentrate in the rest of the paper on value functions that are affine in the state variable. More precisely, we look for affine value functions of the form
\begin{equation}\label{value_affine}
V^{i}(P)=w_i+\int_{\Omega}v_i(\bs{x})P(\bs{x})\,d\bs{x},
\end{equation}
for some unknowns $w_i\in\mathbb{R}$ and $v_i\in\mathbb{X}$.

We observe that the Fr\'{e}chet derivative of the affine functional $V^{i}(P)$ in (\ref{value_affine}) is the linear operator defined by
\begin{equation}\label{Derivative_value}
\big\langle \nabla V^{i}(P),h\big\rangle=\int_\Omega v_i(\bs{x})h(\bs{x})\,d\bs{x},\quad h\in\mathbb{X}.
\end{equation}
Using (\ref{Derivative_value}) and integrating by parts we have
\begin{equation}\label{Derivative_on A}
\big\langle \nabla V^{i}(P), \bs{A}P\big\rangle
=\int_\Omega v_i\bs{A}Pd\bs{x}=\int_\Omega \bs{A^*}v_iPd\bs{x},
\end{equation}
where $\bs{A}^*$ denotes the adjoint of $\bs{A}$ defined for  $v_i\in D(\bs{A}^*)=D(\bs{A})$ by
$$\bs{A^*}v_i= \nabla\cdot(k\nabla v_i)  - \bs{b}\cdot\nabla v_i - c v_i.$$
Using again (\ref{Derivative_value}) we have
\begin{equation}\label{Derivative_on_control}
\big\langle \nabla V^{i}(P),u_j\mathbbm{1}_{\Omega_j}\big\rangle =\int_\Omega v_iu_j\mathbbm{1}_{\Omega_j}d\bs{x}, \quad j=1,\dots,J.
\end{equation}
Finally, the derivatives with respect to $u_i$ of the functionals $\mathcal{G}^{i}(P,u_i)$ and
$$\mathcal{F}^{i}(u_i):=\int_{\Omega}v_iu_i\mathbbm{1}_{\Omega_i}d\bs{x}$$
are given by
\begin{align}
&\big\langle\nabla_{u_i}\mathcal{G}^{i}(P,u_i),h\big\rangle=
\int_{\Omega}\frac{1}{u_i}h\mathbbm{1}_{\Omega_i}d\bs{x},\quad h\in\mathbb{X},\label{Derivative_G}\\
&\big\langle \nabla_{u_i}\mathcal{F}^{i}(u_i),h\big\rangle=
\int_{\Omega}v_ih\mathbbm{1}_{\Omega_i}d\bs{x}, \quad h\in\mathbb{X}.\label{Derivative_F}
\end{align}

We can write now the first-order condition for the maximization of the right hand side in (\ref{HJB}) using (\ref{Derivative_value})-(\ref{Derivative_F}). We have
\begin{equation}\label{first_order}
\int_{\Omega}\frac{1}{u_i}h\mathbbm{1}_{\Omega_i}d\bs{x}
+
\int_{\Omega}v_ih\mathbbm{1}_{\Omega_i}d\bs{x}=0,\quad \forall h\in\mathbb{X}=L^2(\Omega).
\end{equation}
So that
\begin{equation}\label{u_i_v_i}
u_i=-\frac{1}{v_i}\mathbbm{1}_{\Omega_i},\quad i=1,\dots J.
\end{equation}

Substituting the guess (\ref{value_affine}) in the HJB equation
(\ref{HJB}) and using the previous results we have
\begin{equation}\label{igualdad}
\begin{split}
\rho w_i +\rho &\int_{\Omega}v_iP\,d\bs{x}=\\
&\int_{\Omega}\bigl(\log\frac{1}{v_i}-\varphi_iP\bigr)\mathbbm{1}_{\Omega_i}d\bs{x}
+\int_{\Omega}\bs{A^*}v_iPd\bs{x}-\int_{\Omega}\sum_{j=1}^J\frac{v_i}{v_j}\mathbbm{1}_{\Omega_j}d\bs{x}.
\end{split}
\end{equation}

From (\ref{igualdad}) it is clear that function $v_i$ satisfies
\begin{equation}\label{equation_vi}
\bs{A^*}v_i=\rho v_i+\varphi_i\mathbbm{1}_{\Omega_i},\quad i=1\dots, J.
\end{equation}
The scalar $w_i$ can be computed from
\begin{equation}\label{w_i}
w_i=\frac{1}{\rho}\int_{\Omega}\log\frac{1}{v_i}\mathbbm{1}_{\Omega_i}d\bs{x}
-\frac{1}{\rho}\int_{\Omega}\sum_{j=1}^J\frac{v_i}{v_j}\mathbbm{1}_{\Omega_j}d\bs{x}.
\end{equation}

Summarizing we have proved the following proposition
\begin{proposition}
There exists a stationary Markov Perfect Nash Equilibrium of the differential game (\ref{objective})- (\ref{state_equation}), such that the value function of player $i$ has the form (\ref{value_affine}) where $v_i$ is the unique solution of the elliptic problem
\begin{equation}\label{v_i}
\begin{aligned}
\nabla\cdot(k\nabla v_i)  - \bs{b}\cdot\nabla v_i - c v_i-\rho v_i
&=\varphi_i\mathbbm{1}_{\Omega_i},\quad \text{in } \Omega\\
\alpha v_i+ k\nabla v_i\cdot\bs{n}&= {0},\quad\quad\,\,\,\,\text{on }
\partial\Omega
\end{aligned}
\end{equation}
and $w_i\in\mathbb{R}$ is given by (\ref{w_i}). The strategy of  player $i$, $u_i=\varphi_i(P)$ is positive and can be computed by (\ref{u_i_v_i}).

Furthermore, the stationary steady state of the pollution stock is given by the solution of the elliptic problem
\begin{equation}\label{steady_pollution}
\begin{aligned}
\nabla\cdot(k\nabla P)  + \bs{b}\cdot\nabla P - c P
+\sum_{j=1}^Ju_j\mathbbm{1}_{\Omega_j}&=0,\quad \text{in } \Omega\\
\alpha P+ k\nabla P\cdot\bs{n}&= {0},\quad \text{on }
\partial\Omega.
\end{aligned}
\end{equation}
\end{proposition}

\section{Numerical examples}

In this section we present some examples in order to illustrate the effects of the spatial configurations on the strategic behaviour of the agents. We remark that even with this simple formulation of the problem it is evident that on the one hand, the model is able to capture the differences that can be expected in the equilibrium policies. On the other hand, the examples show that, as expected,  the optimal location of emissions as well as its size   depend on the geographical position among the players.

We fixed the following values of the parameters: $c = 0.5$, $\varphi_i=1$, $\rho=0.01$, $k=1$. Also and without any loss of generality we choose $P_b=0$. We remark that with this choice of the parameters all the agents are symmetric except, perhaps, for the geographic relative position. It should be apparent that the qualitative results do not depend on the particular values chosen if we restrict our study to constant coefficients (isotropic diffusion). The case of anisotropic diffusion, although interesting, is out of the scope of this paper.

\begin{example}

In this example we consider a rectangular region $\Omega$ subdivided in two identical countries $\Omega_i$, $i=1,2$.
More explicitly $\Omega=\Omega_1\cup\Omega_2$ with $\Omega_1=[0,0.5]\times[0,1]$ and $\Omega_2=[0.5,1]\times[0,1]$, respectively. We consider that region $\Omega$ is completely isolated from the exterior, that is we put $\alpha=0$ in (\ref{state_equation}). We consider $\bs{b}=\bs{0}$ (no convection) in this example. The left picture in Figure~1 represents the emissions and the right picture represents the stock of pollution at steady state. We fix this convention for the rest of examples in this paper. We can observe that both countries behave in a symmetric manner as corresponds to the completely symmetric geometry.  On the one hand, each country chooses the emission rate symmetrically with respect to the horizontal axis. On the other hand, both countries emit more near the common boundary between $\Omega_1$ and $\Omega_2$. The emission rate decreases as the distance from the common boundary increases. Because both regions are isolated from outside both strategically decide to reduce as much as possible the emission rates at the points in space where the exchange of pollution stock with their neighbours is more difficult. The steady-state levels of the pollution stocks compare as the equilibrium emission rates.

\begin{center}
\begin{figure}[h]
\hskip 1cm
\begin{minipage}[h]{6.5cm}
\mbox{}
\pdfimage height 6cm {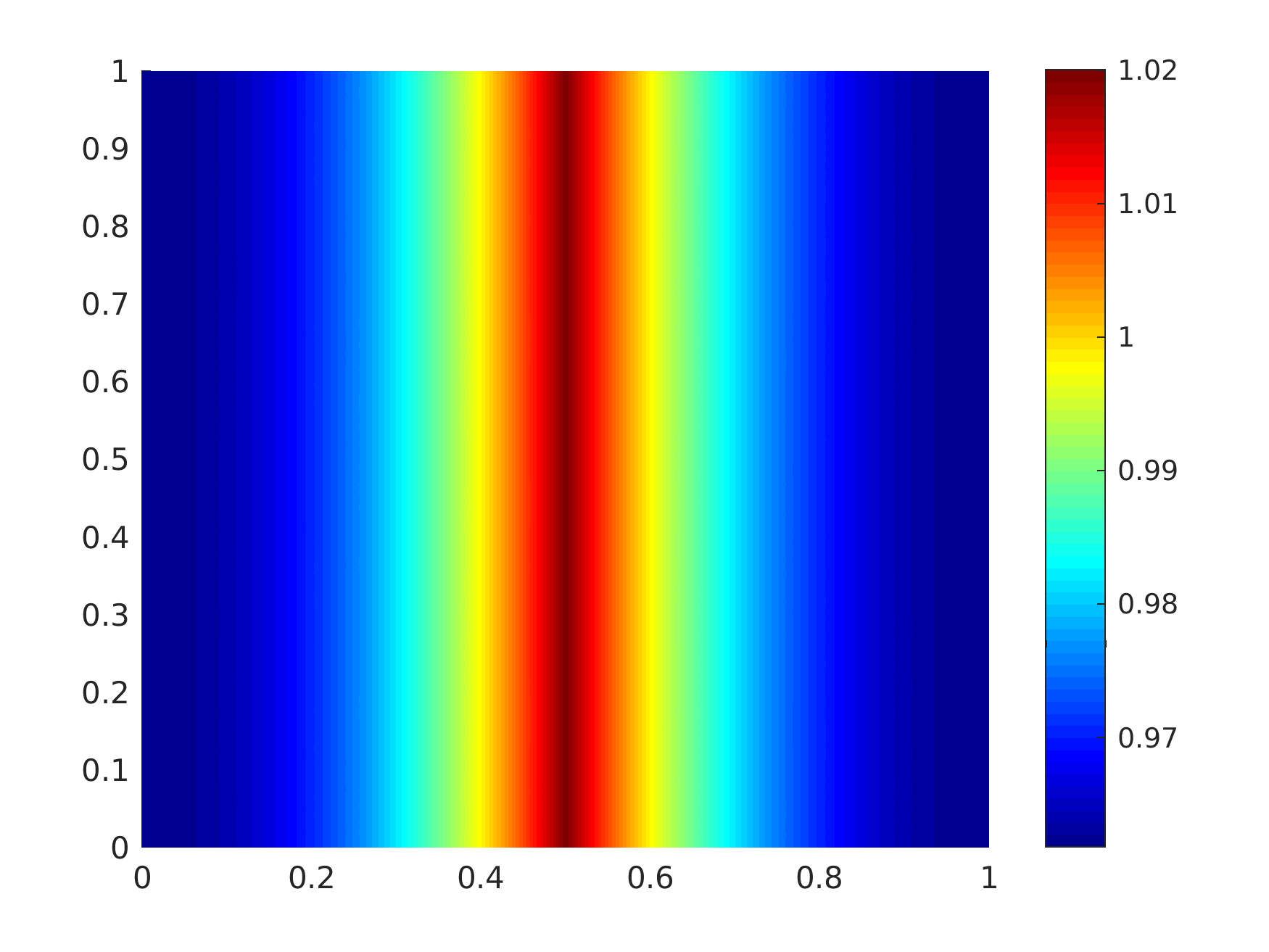}\hfill\mbox{}
\centerline{{\bf Emissions}}
\end{minipage}\hskip 1cm
\begin{minipage}[h]{6.5cm}
\mbox{}
\pdfimage height 6cm {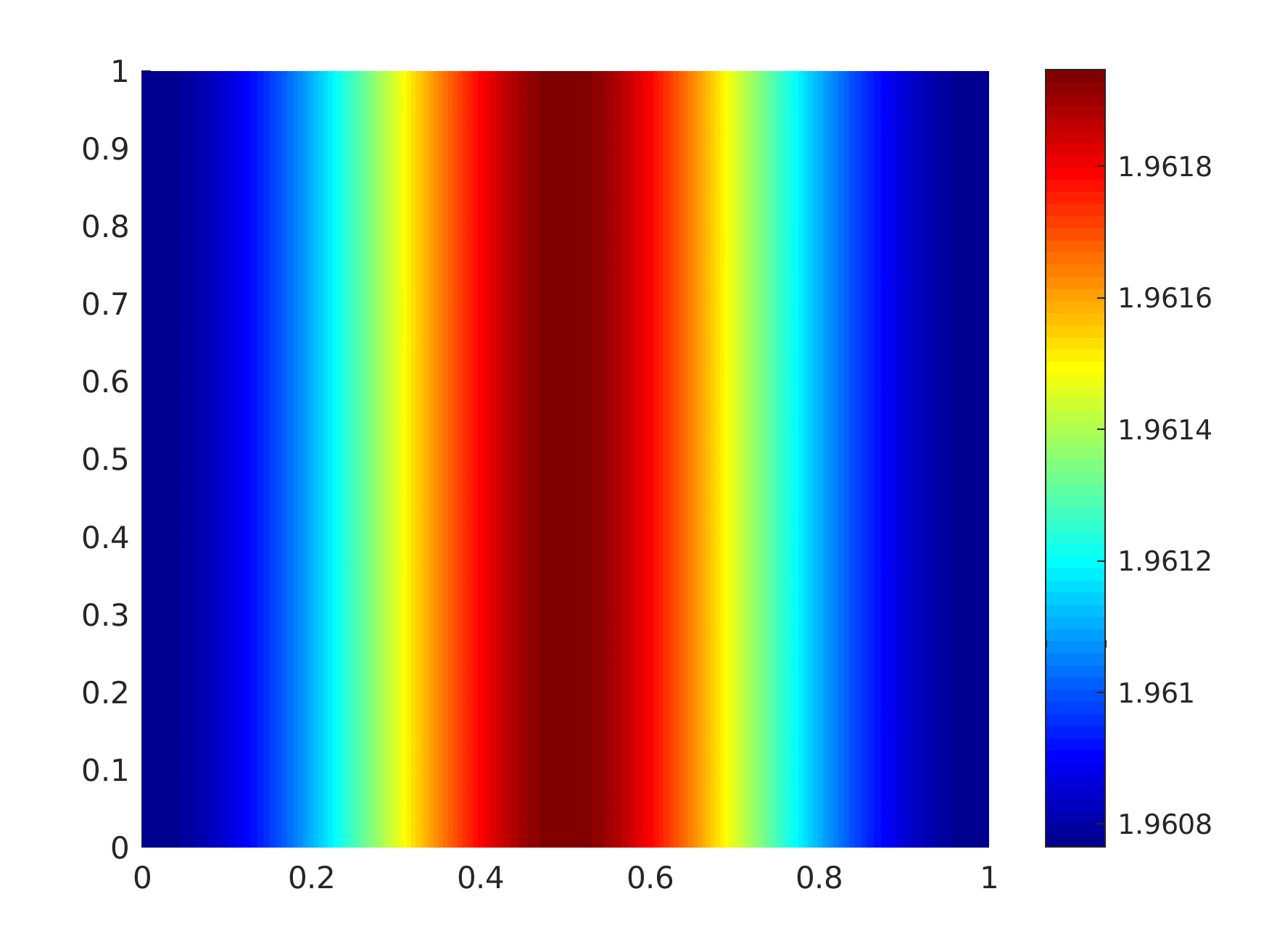}\hfill\mbox{}
\centerline{{\bf Pollution Stock}}
\end{minipage}
\caption{Two isolated symmetric countries. }
\end{figure}
\end{center}

\end{example}

\begin{example}

The geometry of this example is identical to the previous one with the difference that we consider that $\alpha=1$ in the part of the boundary defined by $x=0$ and $x=1$ and $\alpha=0$ on the lines defined by $y=0$ and $y=1$. This choice models a situation in which the region is isolated from the outside in the top and bottom boundaries defined by $y=0$ and $y=1$, and can freely exchange pollution with the clean ($P_b=0$) exterior through the vertical lines $x=0$ and $x=1$. The behaviour is similar to the previous example except for the effect of the open vertical boundaries. The two countries are able to detect that they can get rid of part of the pollution stock at no cost and, consequently, they increase the emission near this open boundary, see Figure~2. Note that the size of the emissions is larger than in Example 1 whereas the stock of pollution is  smaller, clearly reflecting the effect of the open boundaries.

\begin{center}
\begin{figure}[h]
\hskip 1cm
\begin{minipage}[h]{6.5cm}
\mbox{}
\pdfimage height 6cm {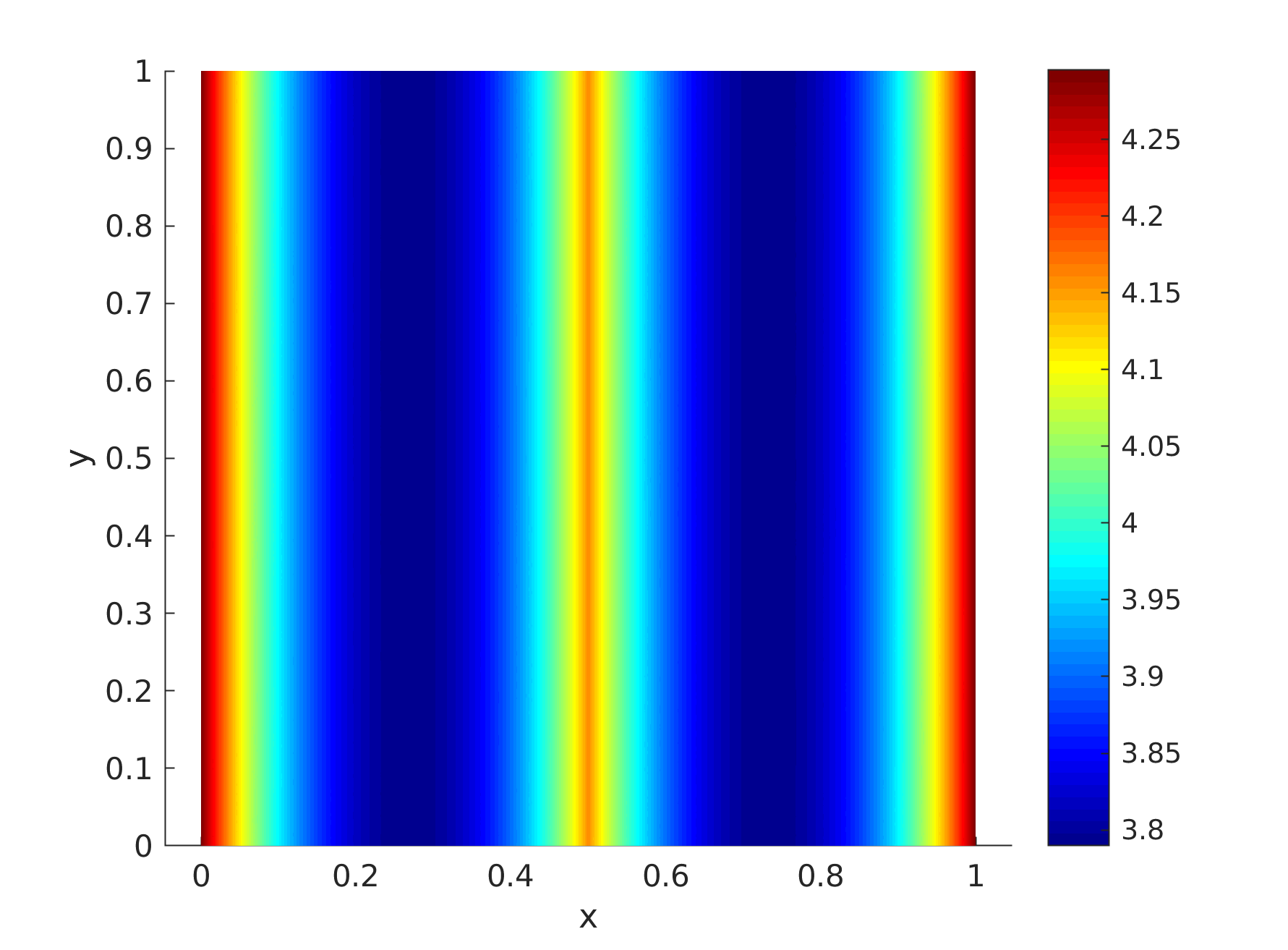}\hfill\mbox{}
\centerline{{\bf Emissions}}
\end{minipage}\hskip 1cm
\begin{minipage}[h]{6.5cm}
\mbox{}
\pdfimage height 6cm {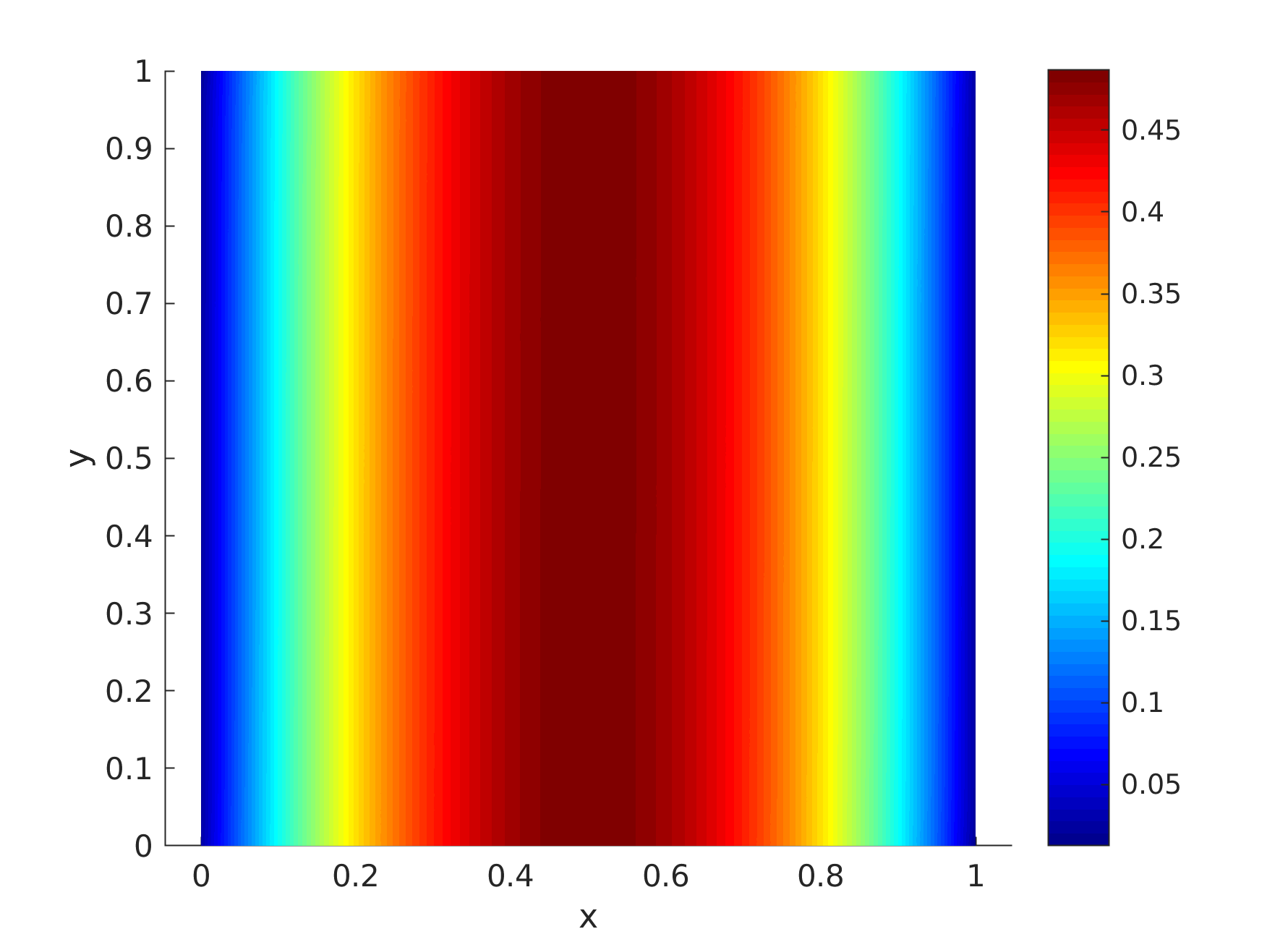}\hfill\mbox{}
\centerline{{\bf Pollution Stock}}
\end{minipage}
\caption{Two symmetric countries. Open boundaries }
\end{figure}
\end{center}
\end{example}

\begin{example}

In this example we consider again a rectangular region $\Omega=[0,2]\times[0,1]$, but this time subdivided in four identical countries with different relative positions. More explicitly the countries are defined by $\Omega_1=[0,0.5]\times[0,1]$, $\Omega_2=[0.5,1]\times[0,1]$, $\Omega_3=[1,2]\times[0,0.5]$ and $\Omega_4=[1,2]\times[0.5,1]$.
 As in the two previous examples $\Omega$ is isolated from the exterior ($\alpha=0$ in (\ref{state_equation})) and there is no convection ($\bs{b}=\bs{0}$). We can observe in Figure~3 that the smaller the number of neighbouring countries, the lower the equilibrium emissions level. This result is consistent with the two  previous examples: Each country knows that its own emissions are being diffused away so they increase the emissions if they have a larger number of neighbours.  This comes from the fact that in this model, the positive effect of the emissions is not shared with other countries, whereas the negative effect of the concentration of pollution is shared through the diffusive state equations. We can also observe  the effect of the relative position on the equilibrium strategies. The emissions in $\Omega_3$ and $\Omega_4$, which are completely symmetric, are much higher than  expected in the proximity of $\Omega_2$, because both countries recognize the increase of emissions in $\Omega_2$. On the contrary, the emissions along the common boundary of $\Omega_3$ and $\Omega_34$ decrease as the distance to $\Omega_2$ increases. As in  Example~1, the steady-state levels of the pollution stocks compare as the equilibrium emission rates.

\begin{figure}[h]
\hskip 1cm
\begin{minipage}[h]{6.5cm}
\mbox{}
\pdfimage height 6cm {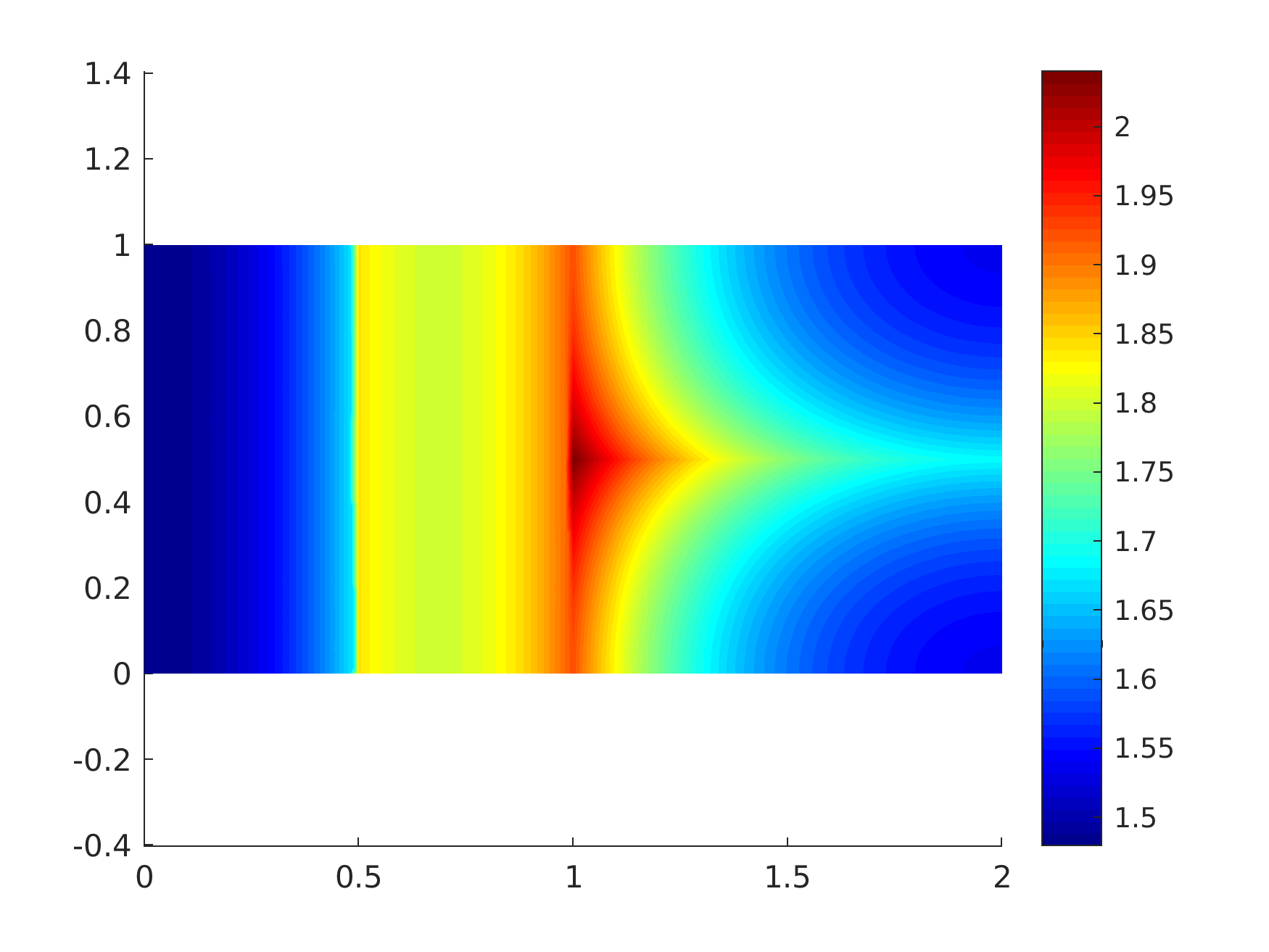}\hfill\mbox{}
\centerline{{\bf Emissions}}
\end{minipage}\hskip 1cm
\begin{minipage}[h]{6.5cm}
\mbox{}
\pdfimage height 6cm {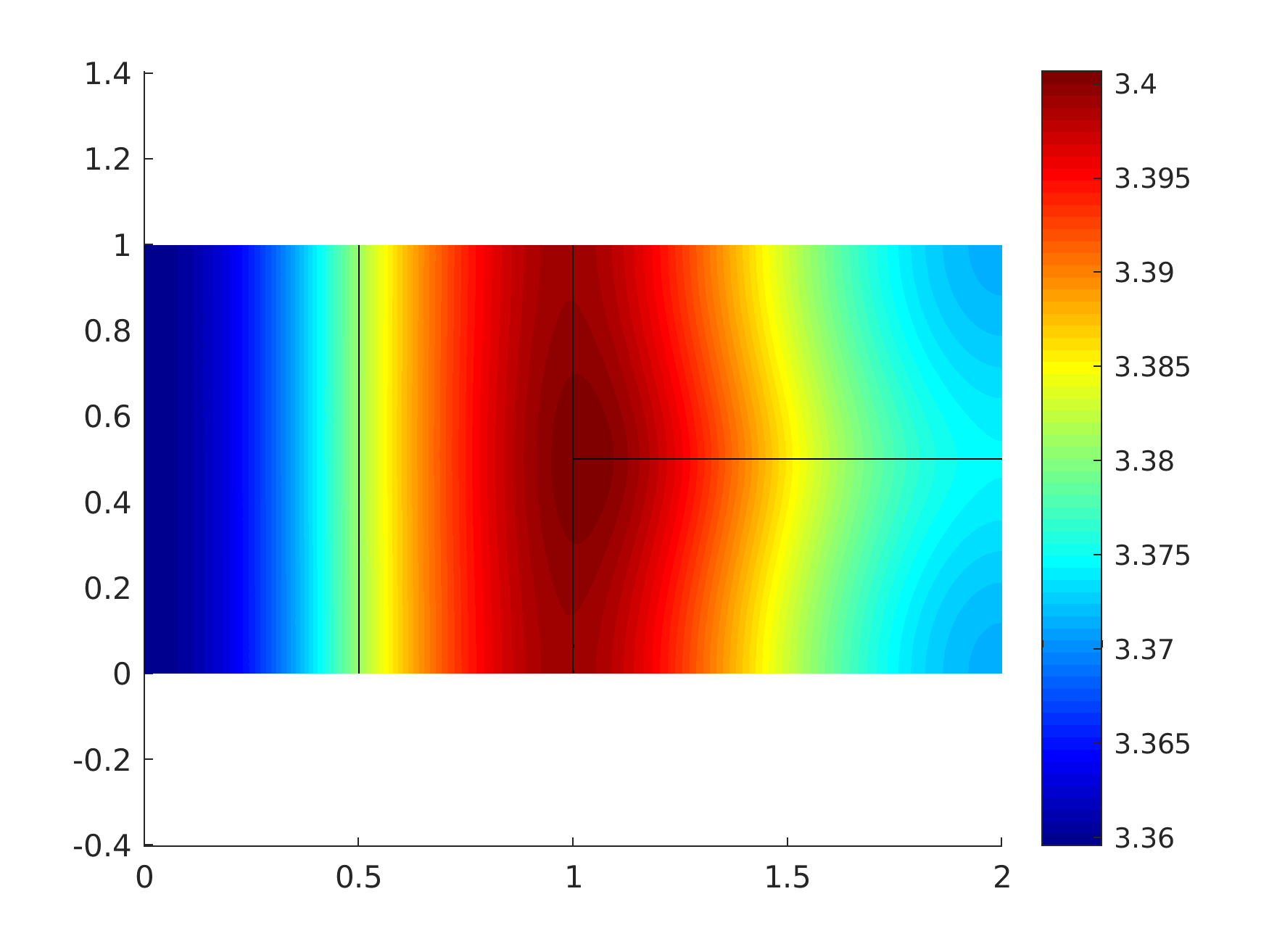}\hfill\mbox{}
\centerline{{\bf Pollution Stock}}
\end{minipage}

\caption{Four countries. Different number of neighbours}
\end{figure}
\end{example}

\begin{example}

Next three examples models a situation in which six, otherwise identical countries, are consecutively positioned along a channel of, say, groundwater forming a (double) L-shaped geometry. The geometry is defined by $\Omega_1=[0,0.5]\times[0,0.5]$, $\Omega_2=[0.5,1]\times[0,0.5]$, $\Omega_3=[0.5,1]\times[0.5,1]$, $\Omega_4=[0.5,1]\times[1,1.5]$, $\Omega_5=[0.5,1]\times[1.5,2]$ and $\Omega_6=[1,1.5]\times[1.5,2]$. In this example  region $\Omega=\bigcup_{j=1}^6\Omega_j$ is isolated from the exterior ($\alpha=0$ in (\ref{state_equation})) and there is no convection ($\bs{b}=\bs{0}$).

\begin{figure}[h]
\hskip 0.5cm
\begin{minipage}[h]{6.5cm}
\mbox{}
\pdfimage height 6cm {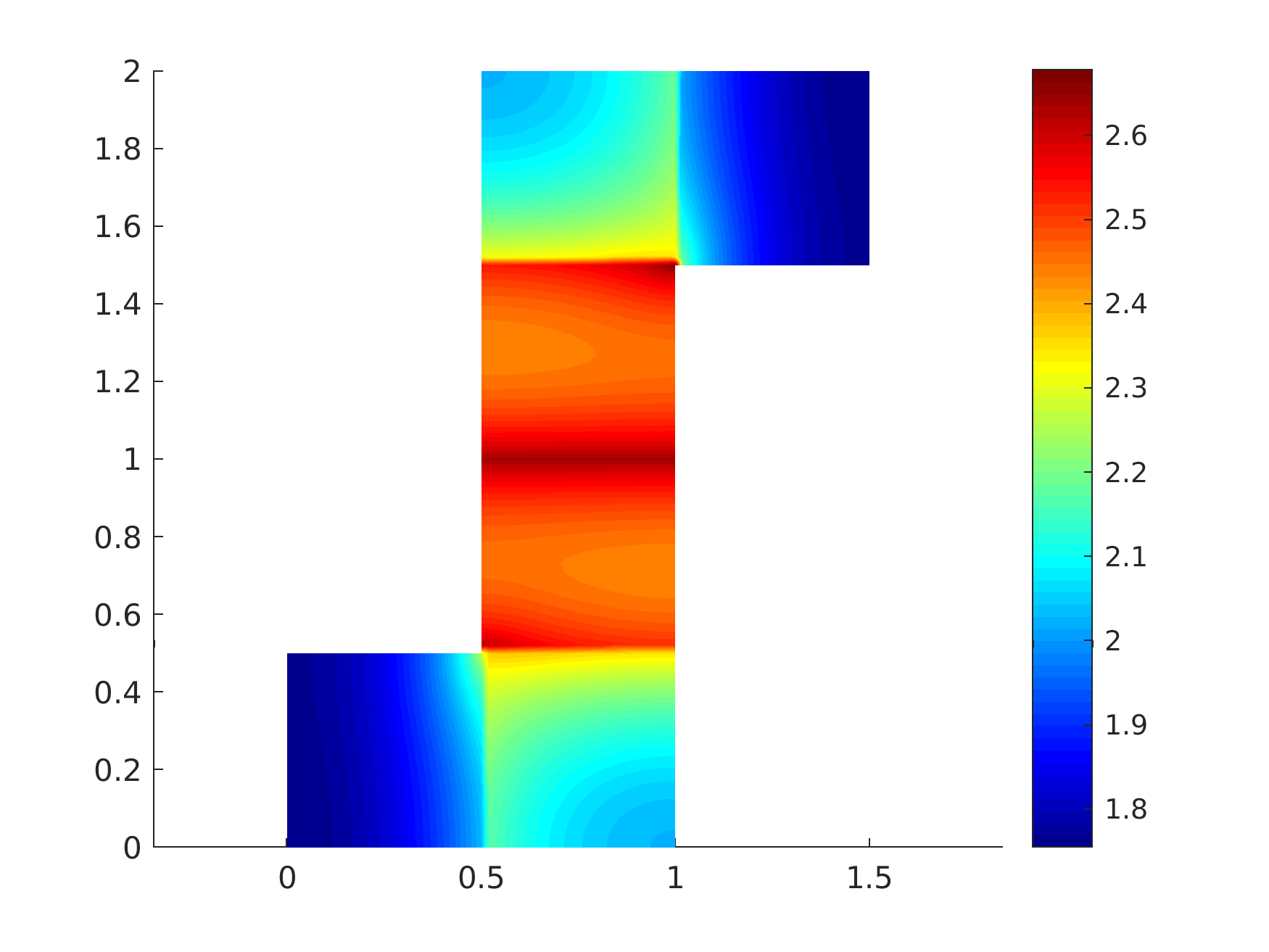}\hfill\mbox{}
\centerline{{\bf Emissions}}
\end{minipage}\hskip 0.75cm
\begin{minipage}[h]{6.5cm}
\mbox{}
\pdfimage height 6cm {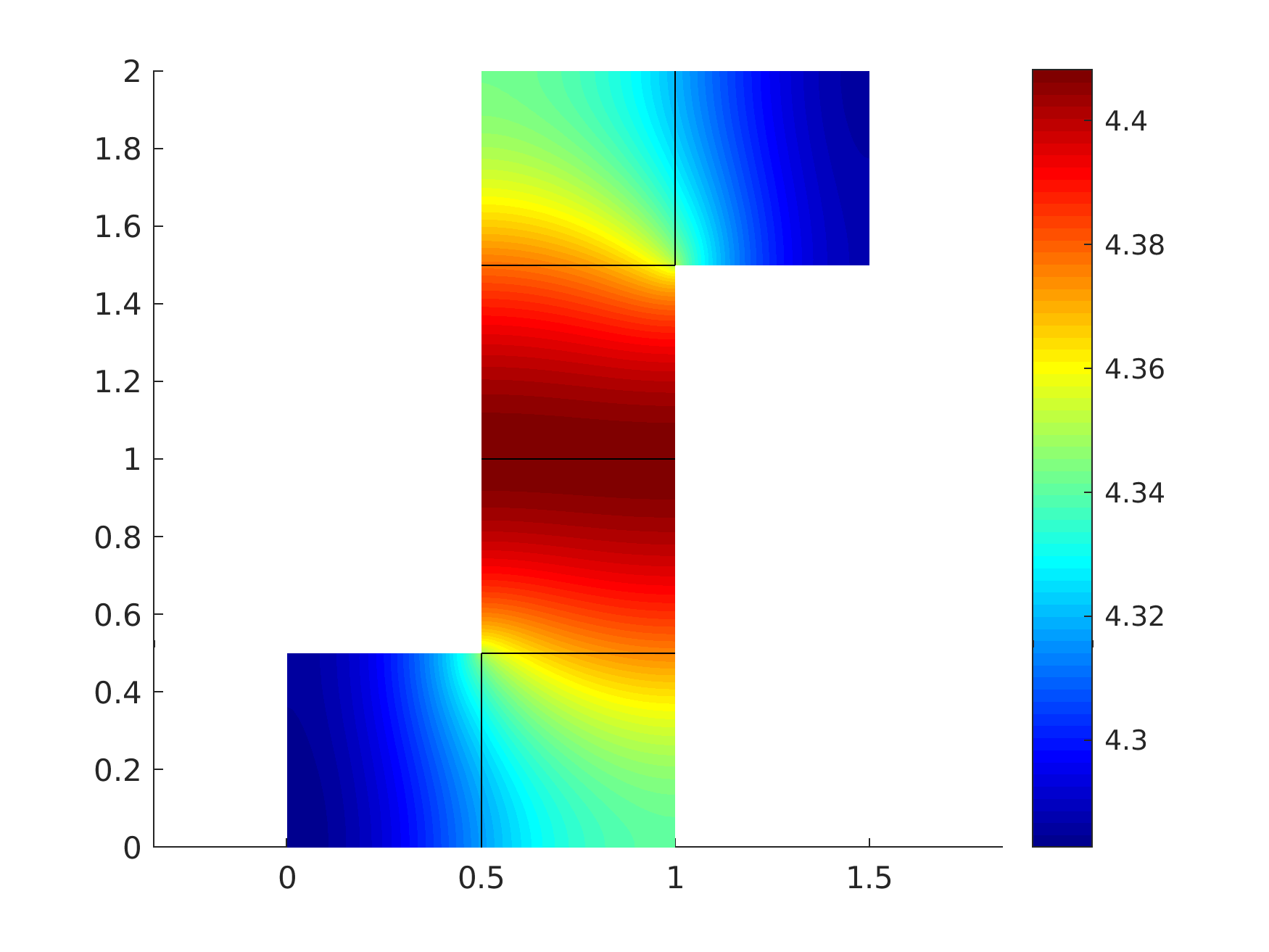}\hfill\mbox{}
\centerline{{\bf Pollution Stock}}
\end{minipage}
\caption{Six countries in a isolated L-shaped domain.}
\end{figure}
\end{example}

In Figure~4 we observe again the dependence of each country emissions on the number of its neighbours. In this picture we can observe a new feature not presented in the previous examples: The dependence on the distance through from the boundary. The emissions in countries $\Omega_3$ and $\Omega_4$ are the highest among the six countries because their neighbours $\Omega_2$ and $\Omega_5$ have themselves another farther neighbour and the diffusion can spread the stock of pollution to a larger and more distant area. We shall use this example to compare with the next and last two examples.

\begin{example}

The geometry of this example is the same as in Example~4 except that now country $\Omega_1$ has an open boundary defined by the line $x=0$. That is, $\alpha=0$ in $\partial\Omega$ except in $\partial\Omega\cap\{x=0\}$ where $\alpha=1$. Figure~5 shows how a small change in the geographical setting can cause a  dramatic change in the equilibrium emissions which is reflected in the change of the steady state of the stock of pollution. Country $\Omega_1$ that has the open boundary takes advantage of this fact emitting at the highest level (see left picture in Figure~5). The emissions of the rest of the countries are lower the longer the distance to the open boundary. However, the stock of pollution follows the opposite pattern: the larger  the distance to the open boundary, the greater  the steady-state stock of pollution.

\begin{figure}[h]
\hskip 0.5cm
\begin{minipage}[h]{6.5cm}
\mbox{}
\pdfimage height 6cm {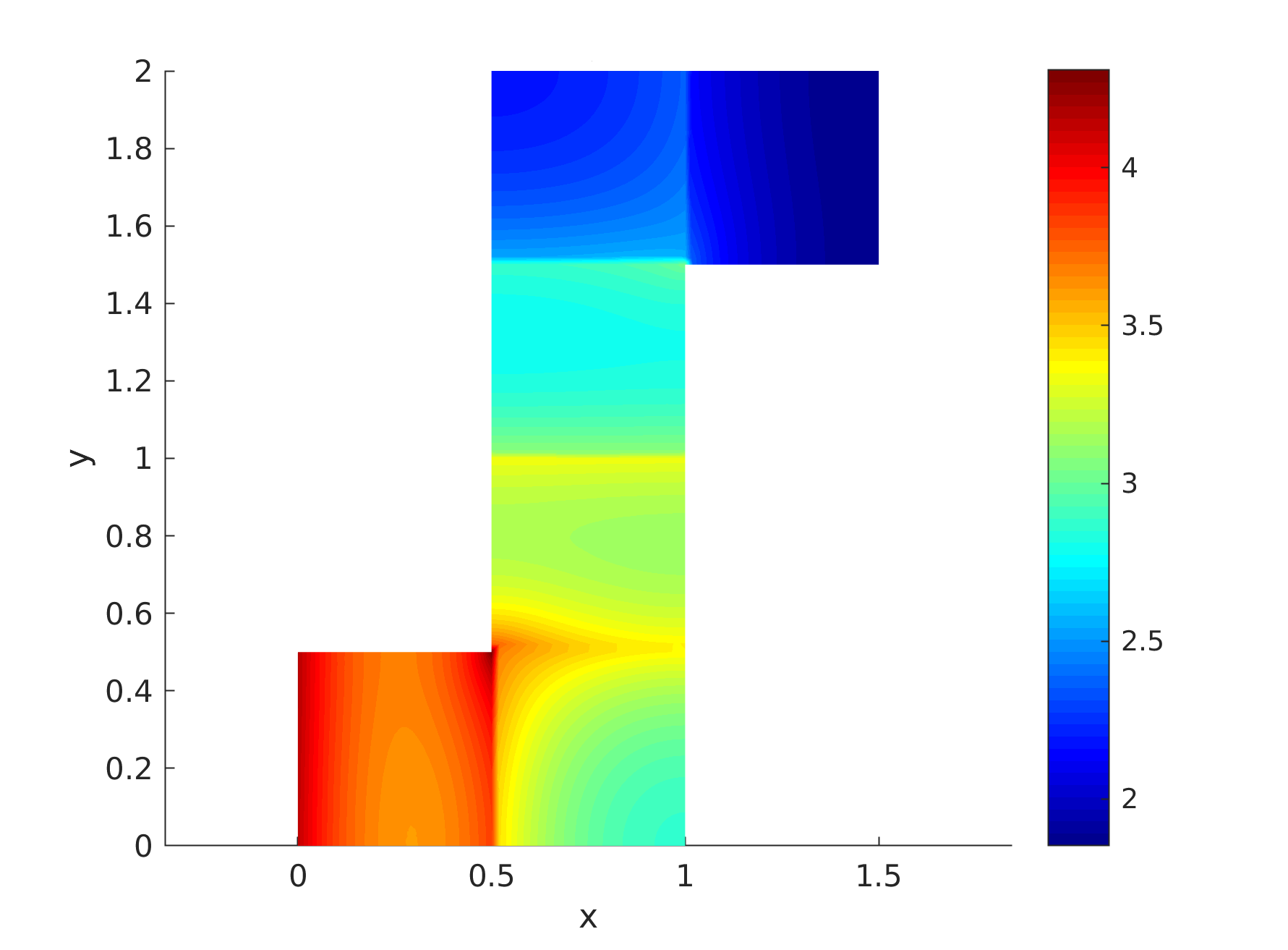}\hfill\mbox{}
\centerline{{\bf Emissions}}
\end{minipage}\hskip 1cm
\begin{minipage}[h]{6.5cm}
\mbox{}
\pdfimage height 6cm {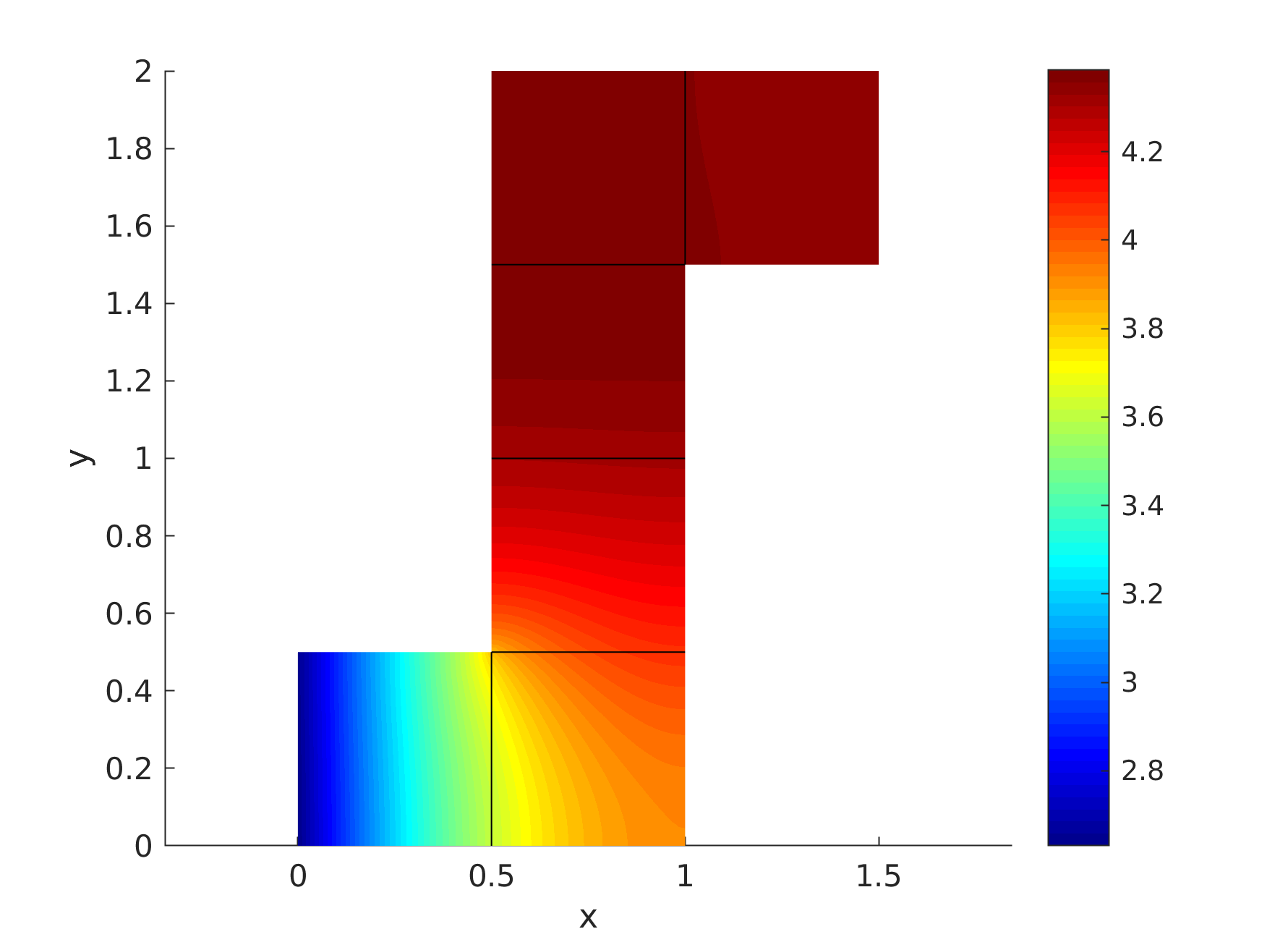}\hfill\mbox{}
\centerline{{\bf Pollution Stock}}
\end{minipage}
\caption{Six countries in an L-shaped domain. Open boundary}
\end{figure}

It is also interesting  to compare the size of emissions in Examples 4 and 5. The following facts are apparent from Figure~4 and Figure~5. First, the emissions are larger in average if $\Omega$ has some part of the boundary open to exchange pollution with the exterior. Of course for this to be true, we have to consider that the exterior is clean ($P_b=0$). However, the size of the steady-state stock of pollution is about the same size in  Figure~4 (right picture) and Figure~5 (right picture). Then although the effect of the open boundary is positive for the overall region $\Omega$, because higher emissions with the same stock of pollution lead to  higher welfare, this is not the case if we individually analyze each one of the  $\Omega_i$.  Clearly, $\Omega_6$ suffers the consequences of the relocation of the stock of pollution which can give, depending on the parameters value, a lower welfare for this country.

\end{example}

\begin{example}

In the last example we consider again the L-shaped domain of Example~5. However, this time we consider a convective flow given by
\begin{equation*}
\bs{b}(\bs{x})=\begin{cases}
(4,0),\quad \bs{x}\in\Omega_1\cup\Omega_6\cup\bigl (\Omega_2\cap\{(x,y),y<x-1\}\bigr)\cup\bigl(\Omega_5\cap\{(x,y),y\ge5/2-x\}\bigr )\\
(0,4),\quad \bs{x}\in\Omega_3\cup\Omega_4\cup\bigl( \Omega_2\cap\{(x,y),y\ge x-1\}\bigr)\cup\bigl(\Omega_5\cap\{(x,y),y<5/2-x\}\bigr).
\end{cases}
\end{equation*}
This means that in this example we are considering an external convective field in the direction of $-\bs{b(x)}$ which is represented in Figure~6. The rest of the parameters and boundary conditions are the same as in Example~5.
\begin{figure}[h]
\hskip 0.5cm
\begin{minipage}[h]{6.5cm}
\mbox{}
\pdfimage height 6cm {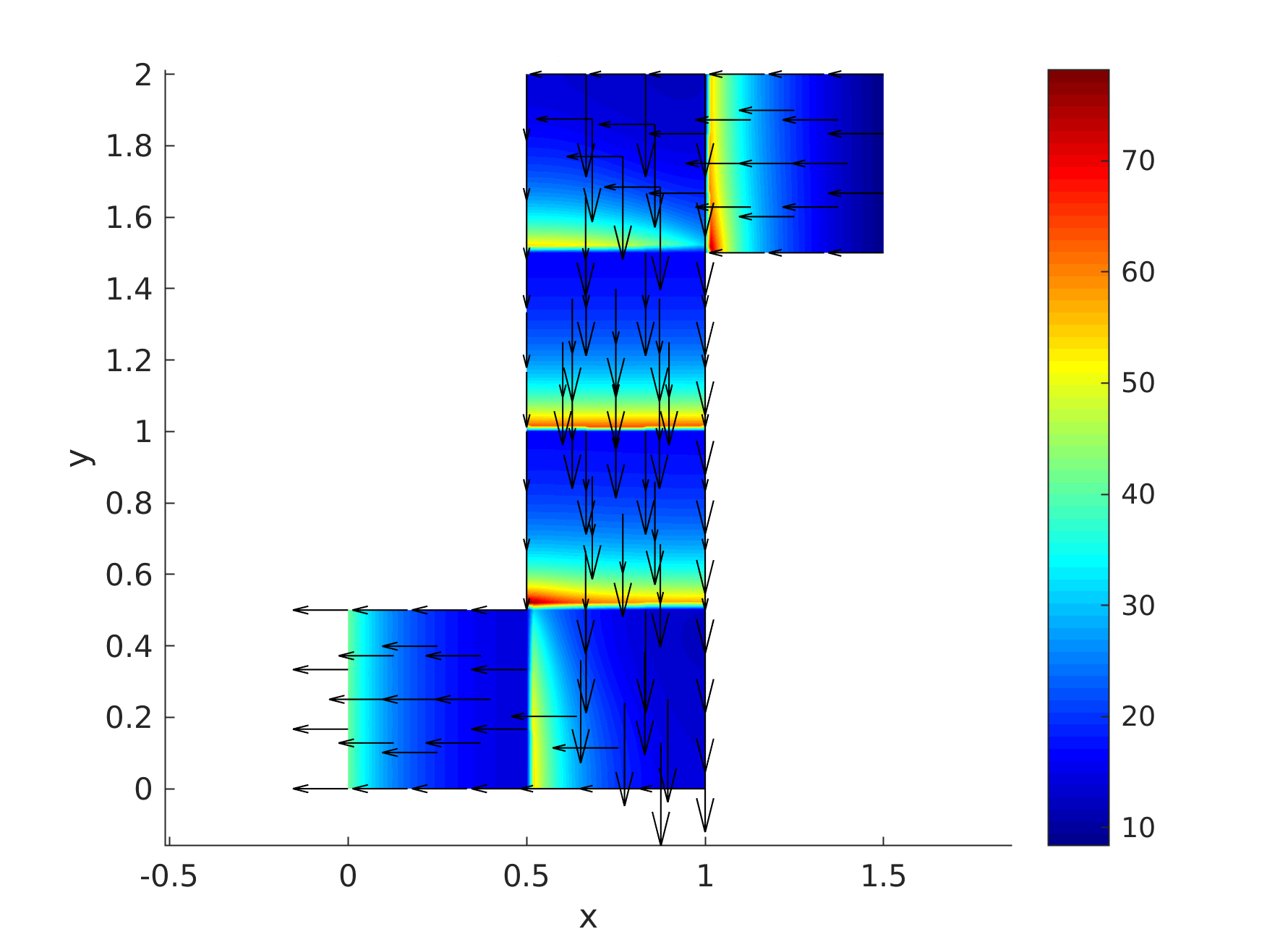}\hfill\mbox{}
\centerline{{\bf Emissions}}
\end{minipage}\hskip 1cm
\begin{minipage}[h]{6.5cm}
\mbox{}
\pdfimage height 6cm {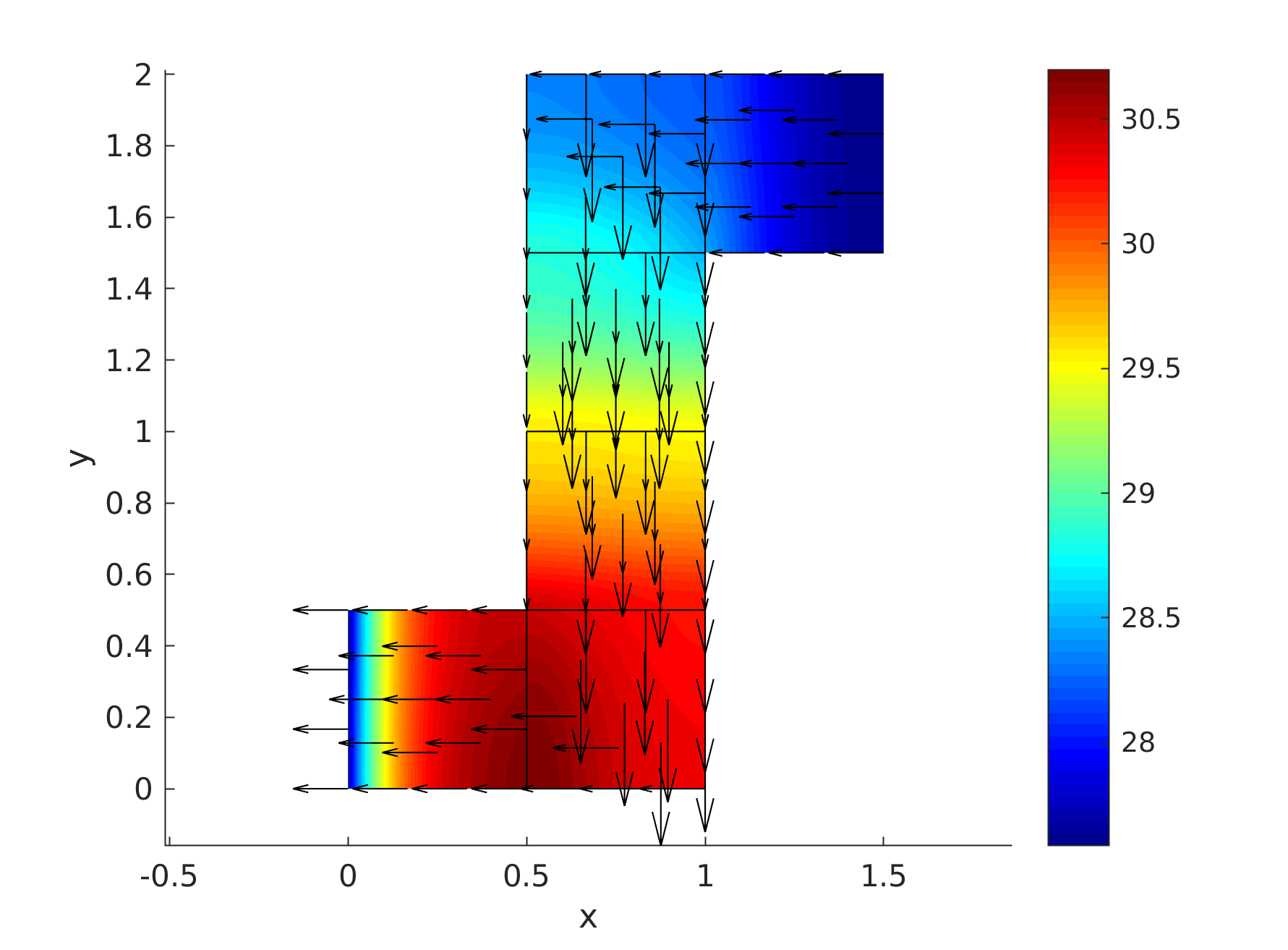}\hfill\mbox{}
\centerline{{\bf Pollution Stock}}
\end{minipage}
\caption{Six countries in an L-shaped domain with convection.}
\end{figure}

We note that the convective field $\bs{b}$ is not parallel to the boundary in $\Omega_1\cap\{(x,y),x=0\}$, where an outward flow is prescribed and in $\Omega_6\cap\{(x,y),x=1.5\}$, where  an inward flow is prescribed. The adjoint problem to be solved in each subdomain $\Omega_i$, $i=1,\dots, 6$, in this particular problem is

\begin{equation*}
\begin{aligned}
\nabla\cdot(k\nabla v_i)  - \bs{b}\cdot\nabla v_i - c v_i-\rho v_i
&=\varphi_i\mathbbm{1}_{\Omega_i},\quad \text{in } \Omega\\
\nabla v_i\cdot\bs{n}-(\bs{b}\cdot\bs{n}) v_i&= {0},\quad\quad\,\,\,\,\text{on }
\Gamma_N\\
\nabla v_i\cdot\bs{n}+(1-\bs{b}\cdot\bs{n}) v_i&= {0},\quad\quad\,\,\,\,\text{on }
\Gamma_R
\end{aligned}
\end{equation*}
where $\Gamma_R=\partial\Omega\cap\{(x,y),x=0\}$ and $\Gamma_N=\partial\Omega\backslash \Gamma_R$.

In Figure~6 we observe that there is a change in scale due to the size of the convective field $\bs{b}$. Qualitatively the effect of the convective flow is clear. The countries concentrate their emissions near the boundary of the downstream neighbour. That is, the countries profit from the fact that pollution is transported away with flow whereas the positive effect of the emissions in the welfare remains in the country. The steady state of the stock of pollution also presents  some interesting and new characteristics. We can clearly see in the right picture of Figure~6 that most part of the stock of pollution is concentrated downstream, even if the size of emissions is similar in the six countries. Interestingly, we also observe  that near the boundary $\Gamma_R$ where the outward flow drains the pollution out of $\Omega$ the stock of pollution is remarkably low due to the open boundary.
\end{example}

\section{Concluding remarks}

This paper studies a transboundary pollution dynamic game that takes into account both the temporal and the spatial dimension of the environmental-economic problem. The spatial-temporal evolution of the stock of a pollutant is described by a diffusion partial differential equation and general boundary conditions are assumed. The main difference with respect to most of the recent literature that adds the spatial aspect in the study of different economic and environmental problems is that instead of a single decision maker, several decision makers and strategic interactions between them are considered. As far as we know  only these two  papers De Frutos \& Mart\'\i n-Herr\'an (2019a, 2019b))  have introduced the  strategic interactions between the decision-makers. These papers
follow a spatial discretization approach to characterize the equilibrium emission strategies of transboundary pollution dynamic games with spatial effects. In this paper we depart from the spatial discretization approach and characterize the feedback Nash equilibrium emission strategies of the $J$-player model formulated in continuous space and continuous time  with two spatial dimensions and one temporal dimension. For a linear-state specification of transboundary pollution dynamic game inspired in J{\o}rgensen \& Zaccour (2001) and De Frutos \& Mart\'\i n-Herr\'an (2019b) we can explicitly  compute the players' value functions and hence, the feedback Nash equilibrium emission strategies. These strategies are constant in time but remarkably they are space-dependent. To the best of our knowledge, we are the first in the literature to characterize the optimal
intraregional distribution of the emissions. The present approach allows us, first, to determine the average total emission in each region (as in De Frutos \& Mart\'\i n-Herr\'an (2019a, b) using the space-discretized model); and second, the optimal spatial location of the emissions of each region. Our analytical results and examples with different geographical configurations show that the spatial aspects play an important role in the determination of the equilibrium emission strategies.

\end{document}